 \input amstex

 \NoBlackBoxes

 \documentstyle{amsppt}
\topmatter

 \title {Newton polyhedra and good compactification theorem}
\endtitle

\rightheadtext{Newton polyhedra and good compactification
}
\footnotetext"*"{The work was partially supported by the Canadian Grant No. 156833-17. }

\author{Askold Khovanskii}
\endauthor
\abstract{A new transparent proof of the well known good compactification theorem for the complex torus $(\Bbb C^*)^n$ is presented. This theorem provides  a powerful tool in enumerative geometry for subvarieties in the complex torus.  The paper also contains an algorithm  constructing   a good compactification for a subvariety in $(\Bbb C^*)^n$  explicitly defined by a system of equations. A new theorem on  a torodoidal like compactification  is stated. A transparent proof of this generalization of the good compactification theorem   which is similar to proofs and constructions from this paper will be presented in  a forthcoming publication}
\endabstract

\affil University of Toronto
\endaffil

\endkeywords

\endtopmatter
\document


\bigskip

\bigskip
\subhead {1. Introduction}
\endsubhead
The paper is dedicated to a simple constructive proof of the good compactification theorem for the group $(\Bbb C^*)^n$ and to related elementary geometry of this group.

A few words about the introduction. We briefly talk about  the ring of conditions $\Cal R(T^n)$ for the complex torus $(\Bbb C^*)^n$ in  the introduction only. This ring suggests a version of intersection theory for algebraic cycles in  $(\Bbb C^*)^n$. We  discuss   a role of the good compactification theorem in this theory and explain how Newton polyhedra are  related to the ring $\Cal R(T^n)$. In section 1.4  we  state a new stronger version of the good compactification theorem.  In the sections 1.5  and 1.6 we summarize the remainder of the paper and fix some notation.
\subsubhead {1.1. The ring of conditions}
\endsubsubhead
The good compactification theorem for a spherical homogeneous space $U$  allows to define{\it the ring of conditions of $U$} [3]. In the paper we consider  the case when $U$ is a complex torus $(\Bbb C^*)^n$ equipped with the natural action of the torus on itself.

Following the original ideas of Schubert, in the early 1980s De Concini and Procesi developed an intersection theory for algebraic cycles  in a symmetric homogeneous space  [3].
Their  theory, named the {\it ring of conditions of $U$}, can be automatically generalized to a spherical homogeneous space $U$. De Concini and Procesi showed that the description of such a ring can be reduced to homology rings (or to Chow rings) of an increasing chain of smooth  $G$-equivariant compactifications of~$U$.

The ring of conditions $\Cal R(U)=\Cal R_0(U)+\dots +\Cal R_n(U)$ is a commutative graded ring with homogenious components of degrees $0,\dots, n$ where $n$ equals to the dimension of $U$. The component $\Cal R_k(U)$ consists of algebraic cycles $Z_k$ in $U$  of codimension $k$, considered up to an equivalents relation $\sim$. An element in $Z_k$ is a formal sum  of  algebraic subvarieties of codimension $k$ taken with with integral coefficients.

1) Two algebraic cycles $X_1,X_2\in Z_k$ are {\it equivalent} $X_1\sim X_2$ (and define the same element in $\Cal R_k(U)$) if for any cycle $Y\in Z_{n-k}$ for almost any  $g\in G$ the  cycles $X_1$ and $gY,$ as well as cycles  $X_2$ and $gY,$ have finite number of points of intersections and these numbers computed with appropriate multiplicities are equal.

2) Assume that:  a) $X_1, X_2\in Z_k$  and $X_1\sim X_2$; b) $Y_1, Y_2\in Z_m$ and $Y_1\sim Y_2$. For almost any element $g\in G$ for $i=1,2$ let $W_{i, g}$  be the cycle obtained by taken with appropriate coefficients components of codimension $k+m$ in the  intersection of  $X_i$ and  $g Y_i$. Then the cycles $W_{i, g}$ and $W_{i, g}$ are equivalent  and they define the same element in $\Cal R_{k+m}(U)$.

The addition $+$ in the ring $\Cal R(U)$ is induced by the addition in the group of cycles. The multiplication $*$  in the ring $\Cal R(U)$ is induced by intersection of cycles. The multiplication in $\Cal R(U)$ is well-defined because of the property 2).

The component $\Cal R_0(U)$  is isomorphic to $\Bbb Z$, each element in $\Cal R_0(U)$ is the variety $U$ multiplied by an integral coefficient.

The component $\Cal R_n(U)$ is isomorphic to $\Bbb Z$. Its elements are linear combinations of points in $U$ with the following equivalence relation: $\sum k_i a_i\sim \sum m_j b_j$ if  $\sum k_i = \sum m_j$.

The pairing $F_k:\Cal R_k\times \Cal R_{n-k}\rightarrow \Bbb Z$ can be defined using the multiplication in the ring $\Cal R(U)$ and  the isomorphism $\Cal R_n\sim \Bbb Z$.  This pairing  is non degenerate, i.e. for any $a_k\in \Cal R_k\setminus \{0\}$ there is  $b_{n-k}\in \Cal R_{n-k}$ such that $F_k(a_k,b_{n-k})\neq 0$ (the non degeneracy follows automatically from the definition of the ring $\Cal R(U)$).

In subsection 1.2  we  comment on the role of the good compactification theorem for this intersection theory  in the  case when a homogeneous space $U$ is the  complex torus $(\Bbb C^*)^n$ acting naturally on itself.

\subsubhead {1.2. The ring $\Cal R (T^n)$ of conditions of  $T^n=(\Bbb C^*)^n$}
\endsubsubhead
The construction of the ring and all its descriptions are based on good compactification theorem
[4--5], [7--8], [10], [17].

A complete  toric variety $M\supset (\Bbb C^*)^n$ is a {\it good compactification for an  algebraic variety $X\subset (\Bbb C^*)^n$} of codimension $k$ if the closure  of $X$ in $M$ does not intersect any orbit of the $(\Bbb C^*)^n$ action  on  $M$ whose dimension is smaller than $k$. The toric variety $M$ is a good compactification for an algebraic cycle  of codimension $k$ if it is a good compactification of all varieties of codimension $k$ which appear in the cycle with nonzero coefficients.
\proclaim {Good compactification theorem ([3], [5], [17], [18])} One can find  a good compactification  for any given algebraic subvariety $X$ in $(\Bbb C^*)^n$.
\endproclaim
Let $M_1,M_2$ be toric varieties and $\pi:M_1\rightarrow M_2$ be a proper equivariant map. Let $M_2$ be  a good compactification for a variety $X$. Then $M_1$ also is a good compactification for $X$. Thus the theorem implies that one can find  a smooth projective toric variety $M$ which provides a good compactification  for any given algebraic subvariety $X$ in $(\Bbb C^*)^n$.

Assume that a smooth projective toric variety $M$ is a good compactification for cycles $X_1,X_2\in Z_k$. Then $X_1\sim X_2$ if and only if the closure of these cycles in $M$ define equal elements in the homology group  $H_{2(n-k)}(M,\Bbb Z)$. Moreover for any element  $x\in H_{2(n-k)}(M,\Bbb Z)$ one can find a cycle $X\in Z_k$ such that $M$ is a good compactification for $X$ and the closure of $X$ in $M$ represents the element $x$.

Assume that  the variety $M$ is a good compactification for cycles $X\in Z_k$ and $Y\in Z_m$ and the closure of these cycles define elements $x\in H_{2(n-k)}(M,\Bbb Z)$  and $y\in H_{2(n-m)}(M,\Bbb Z)$. Let  $g$ be a generic element in $(\Bbb C^*)^n$. Let $W_g \in Z_{k+m}$ be the cycle   corresponding to the intersection of $X$ and $gY$. Then $M$ is a good compactification for $W_g$. Moreover the closure of $W_g$ in $M$ defines the element $w\in H_{2(n-k-m)}(M,\Bbb Z)$ which is independent of a choice of generic element $g$ and equal to the intersection of cycles $x$ and $y$ in the homology ring $H_{*}(M,\Bbb Z)$.

These statements show that  description ring $\Cal R(T^n)$ can be reduced to the description of the homology rings of smooth projective varieties (and to  description of the behavior of these rings under proper equivariant  maps between toric varieties) (see [3--5], [10]). Related material can be found in [7--8] and [17].

\subsubhead {1.3. Newton polyhedra and the ring $\Cal R(T^n)$}
\endsubsubhead
Let us start with the following striking  connection  of Newton polyhedra  to the ring $\Cal R(T^n)$. Consider two hypersurfaces $\Gamma_1, \Gamma_2\subset (\Bbb C^*)^n$ defined by equations $P_1=0$, $P_2=0,$ where $P_1$ and $P_2$ are Laurent polynomials. Then {\it the cycles $\Gamma_1, \Gamma_2\in Z_1$ represent the same element in $\Cal R_1(T^n)$ if and only if the Newton polyhedra $\Delta(P_1)$ and $\Delta(P_2)$ are equal up to a shift.}

Any element in $\Cal R_1(T^n)$ can be represented as a formal difference of two hypersurfaces. Thus {\it the component $\Cal R_1(T^n)$ can be identified with the group of virtual convex polyhedra with integral vertices}.

{\it The ring $\Cal R(T^n)$ is generated by $\Cal R_0(T^n)\sim \Bbb Z$ and $\Cal R_1(T^n)$ }(since the homology ring of a smooth projective toric variety $M$ is generated by $H_{2n}(M,\Bbb Z)\sim \Bbb Z$ and $H_{2n-2}(M,\Bbb Z)$). Thus it is not surprising that the ring $\Cal R(T^n)$ can be described using only the geometry of convex polyhedra with integral vertices  [4], [10].

Let  $\Gamma_1,\dots,  \Gamma_n\subset (\Bbb C^*)^n$ be hypersurfaces defined by equations $P_1=0,\dots, P_n=0$ where $P_1,\dots, P_n$ are Laurent polynomials with Newton polyhedra $\Delta_1,\dots,\Delta_n$.  Then the intersection number of the cycles $\Gamma_1, \dots, \Gamma_n$ in the ring of conditions $\Cal R(T^N)$ is equal to $n !$ multiplied by the mixed volume of $\Delta_1,\dots,\Delta_n$ (see [10]). This statement is a version of the famous Bernshtein--Koushnirento theorem (see [1--2], [12--13], [16])  also known as BKK (Bernshtein--Koushnirenko--Khovanskii) theorem. The statement shows that BKK theorem applies naturally to the ring $\Cal R(T^n)$.
\subsubhead {1.4. Strong version of the good compactification theorem}
\endsubsubhead
In this section we announce a toroidal like compactification theorem which is a stronger version of the good compactification theorem for complex torus.

Consider a triple $(Y,D,a)$, where $Y$ is a normal $n$-dimensional variety, $D\subset Y$ is a Weil divisor and $a\in D$ is a point. We say that $(Y,D,a)$ is a {\it pointed toroidal triple} if there is an affine toric variety $Y_1\supset (\Bbb C^*)^n$ with a divisor $D_1=Y_1\setminus (Y_1\cap  (\Bbb C^*)^n)$ and with a null orbit $a_1\in D_1$ such that the triple $(Y,D,a)$ is locally analytically equivalent to the triple $(Y_1,D_1,a_1)$ in  neighborhoods of the points $a$ and $a_1$ (about toroidal embedding see [11]).

Let $M\subset(\Bbb C^*)^n$ be a good compactification for $X\subset (\Bbb C^*)^n$, let  $\Cal X$ be the closure of $X$ in  $M$ and let  $\Cal D$ be  $\Cal X\cap D$ where $D =M\setminus (\Bbb C^*)^n$. Consider the normalization $\pi: \tilde \Cal X\rightarrow X$ and let $\tilde D$ be $\pi^{-1}(D)$.

We say that a good compactification $M$ for a subvariety $X$  of codimension $k$ in $ (\Bbb C^*)^n$ gives a {\it toroidal like compactification $\Cal X$  of $X$ if for any point $a\in \Cal X\cap O,$ where   $O\subset M$ is a $k$-dimensional orbit, and for any point $\tilde a\in \tilde \Cal X$ such that $\pi(\tilde a)=a$ the triple $(\tilde \Cal X, \tilde \Cal D,\tilde  a)$ is a pointed toroidal triple.
\proclaim{Toroidal like compactification theorem} For any given algebraic subvariety $X$ in $(\Bbb C^*)^n$ one can find  a toric variety $M\supset (\Bbb C^*)^n$ providing a toroidal like compactification of $X$.
\endproclaim

A transparent proof of this new theorem which is similar to proofs and constructions from this paper will be presented in a forthcoming article.

Consider any toroidal like compactification for a variety $X \subset (\Bbb C^*)^n$. Let $\tilde a\in \tilde \Cal X$ be a point as above and let $\tilde \Cal V\subset \tilde \Cal X$ be a small neighborhood   of $\tilde a$ in $\tilde \Cal X$.  {\it  The image $\Cal V=\pi(\tilde \Cal V)\subset \Cal X$} (which is a piece of a branch of $\Cal X$ passing through the point $\pi(\tilde a)$) {\it can be represent by
converging multidimensional Laurent power series}. From a tropical geometry  point of view domains of convergency of such Laurent series ``cover almost all the infinity of $X$". Such series could be useful for enumerative geometry.

\subsubhead{1.5. Summary of the paper}
 \endsubsubhead
A short proof of the good compactification theorem is presented in the section 7.
We use  {\it convenient compactifications  for a system of equations}
(section 5) and  {\it developed systems of  equations}
 (section 6).  Appropriate {\it regular sequences} (section 2) allow to reduce an arbitrary variety $X$ to a complete intersection $Y$, such that $X\subset y$, $\dim X=\dim Y$ (section 3).

Sections 2 and 4 are not used directly in the proof. In the section 2 we  explain how to reconstruct the dimension of a variety $X$ from the set of Newton polyhedra of all Laurent polynomials in an ideal defining $X$. Our proof of the good compactification theorem is based on similar observations.

The proof heavily uses developed systems of  $k$ equations in $(\Bbb C^*)^n$. In Newton polyhedra theory there are many results on such systems with $k=n$. In section 4 we recall  these results.

An algorithm  constructing   a good compactification for a subvariety  $X$ in $(\Bbb C^*)^n$  explicitly defined by a system of equations is presented in section 10. It is based on elementary results from elimination theory (section 8) which allows to define  by explicitly written system of equations a projection $\pi (X)$  of $X$ on an $(n-1)$-dimensional subtorus (section 9).

In section 12 a modified algorithm is presented. We assume there that the codimension of $X$ is given  and we make arbitrary generic choices in order for the construction to work. This modification is based on the original algorithm and on auxiliary results presented in section 11.

\subsubhead {1.6. Complex torus, its subgroups and factor groups}
\endsubsubhead
Here we fix some notation.

The standard $n$-dimensional complex torus we denote by $(\Bbb C^*)^n$. We denote by $\Cal L_n$ the ring of regular functions on $(\Bbb C^*)^n$  consisting of Laurent polynomials (i.e. of linear combinations of characters of $(\Bbb C^*)^n$). We denote by:

$\Lambda\sim \Bbb Z^n$ -- the lattice of characters;

$M\sim \Bbb R^n$ -- the space  $\Lambda\bigotimes_{\Bbb Z}\Bbb R$ of characters;

$\Lambda^*\sim \Bbb (Z^n)^*$ -- the lattice of one-parameter subgroups;

$N\sim (\Bbb R^n)^*$ -- the space  $\Lambda^*\bigotimes_{\Bbb Z}\Bbb R$ of one parameter subgroups.

\smallskip
A subspace $L\subset M$ is {\it rational} if it is generated by elements from $\Lambda$. With a rational
space $L$ one can associate the following objects:
\smallskip

1) The sublattice $\Lambda (L)$  of characters in $L$ equal to $\Lambda \cap L$.   The sub-ring $\Cal L_L$ of the ring $\Cal L_n$ consisting of Laurent polynomials  whose Newton polyhedra belong to  $L$.
\smallskip

2) The connected subgroup  $H(L)$ in $(\Bbb C^*)^n$ whose Lie algebra is spanned by vectors from $\Lambda^*$ orthogonal to the space $L$. The subgroup  $H(L)$ is a subset in $(\Bbb C^*)^n$ where all characters from $\Lambda (L)$ equal to one. Its dimension $\dim_\Bbb C H(L)$  equals to $n-\dim_\Bbb R L$.

\smallskip
3) The factor-group $F(L) =(\Bbb C^*)^n/H(L)$ and the factorization map $\pi: (\Bbb C^*)^n\rightarrow F(L)$.  The map $\pi^*$ identifies the lattice of characters of  $F_L$ with the sublattice $\Lambda(L)$, and it identifies the ring of regular functions on $F(L)$ with the ring $\Cal L_L$.

\subsubhead {1.7. Acknowledgements}
\endsubsubhead
I would like to thank Boris Kazarnoskii for frequent  long and useful  discussions.  I also am
grateful to Feodor Kogan who edited my English.

\subhead {2. Newton polyhedra and dimension of an algebraic variety}
\endsubhead
Dimension of an algebraic variety $X\subset (\Bbb C^*)^n$ is its key invariant for the good compactification theorem.  Assume that $X$ is defined by an ideal $I$ in the ring of Laurent polynomials.  It is natural to ask the following question: {\it is it possible(at least theoretically)  to determine the dimension of the variety $X$ knowing  the set $\Delta(I)$ of Newton polyhedra $\Delta(P)$ of all Laurent polynomials  $P\in I$?}. Below  we provide a positive  answer to this question.  Later we will not use this result directly, but  our approach to the good compactification theorem is based on simple arguments used in this section.

Let   $X\subset (\Bbb C^*)^n$ be the algebraic variety   defined by an ideal  $I\subset \Cal L_n$.

\proclaim {Lemma 1} If $\dim X=m,$ then for any rational subspace $L\subset M$  such that $\dim_{\Bbb R}L=m+1$ there is a Laurent polynomial $P\in I$ whose Newton polyhedron  $\Delta(P)$ belongs to $L$.
\endproclaim
\demo{Proof} Let $L\subset M$ be an $(m+1)$-dimensional rational space. Then the factor-group $F(L)=(\Bbb C^*)^n/H(L)$ has dimension $(m+1)$. Let $\pi: (\Bbb C^*)^n\rightarrow F(L)$ be the natural factorization map. The image $\pi(X)\subset F(L)$ of $X$ has dimension $\leq m$ and it has to belong to some algebraic hypersurface $Q$ in the $(m+1)$-dimensional group $F(L)$. By definition the function $\pi^*(Q)$ vanishes on $X$.   Thus $\pi^*(Q)$ belongs to the radical of the ideal $I$, i.e. there is a natural number $p$ such that $(\pi^*Q)^p\in I$.  The Newton polyhedron $\Delta$ of $\pi^*(Q)$ belongs to $L$. Now one can choose $P= (\pi^*Q)^p$ since $\Delta (P)=p \Delta  \subset L$.
\enddemo

\proclaim {Lemma 2} If $\dim X=m,$ then for a generic rational subspace $L\subset M,$ with $\dim L=m,$ there is no $P\in I$ such that  $\Delta(P)$ belongs to $L$.
\endproclaim

\demo{Proof} For a generic rational subspace $L\subset M$ with $\dim L=m$ the image of $X$ in the $m$-dimensional factor-group $F(L)$ has dimension $m$ and it can not belong to any algebraic hypersurface in $F(L) $. Thus there is no Laurent polynomial $P\in I$ such that $\Delta(P)\subset L$.
\enddemo

Lemmas 1, 2  allow (at least theoretically) to determine the dimension of $X$ by Newton polyhedra of all Laurent polynomials from an ideal $I$ which defined $X$.

\proclaim{Theorem 3} If $\dim X=m,$ then any rational space $L\subset M$ with $\dim L=m+1$ contains the Newton polyhedron $\Delta(P)$  of some  $P\in I$ and a generic rational space $L\subset M$ with $\dim L=m$ does not contain any such polyhedron.
\endproclaim

\subhead {3. Good compactification and regular sequences}
\endsubhead
One can strongly modify an algebraic variety $X$ preserving its dimension. In this section we will talk abut complete intersections $Y$ containing $X$ such that $\dim Y=\dim X$. We are interested in such varieties $Y$ because of the following obvious observation.
\proclaim {Lemma 4} If $X \subset Y\subset (\Bbb C^*)^n$ and $\dim X=\dim Y$ then any  good compactification for $Y$  is a good compactification for $X$.
\endproclaim

\demo{Proof} Let $M\supset \Bbb (C^*)^n$ be a good compactification for $Y$. Then the closure of $Y$ in $M$ does not intersect any orbit $O\subset M$ such that $\dim O<n-\dim Y=n-\dim X$.  The closure of $X$ in $M$ also does not intersect any such orbit since $X\subset Y$.
\enddemo

Let us recall the following definition.

\definition{Definition 1} A sequence $P_1,\dots,P_k$ of Laurent polynomials on $(\Bbb C^*)^n$ is  {\it regular} if for $i=1,\dots,k$ the following conditions $(i)$ hold.

 (1):   $P_1$ is not identically equal to zero on $(\Bbb C^*)^n$.

$(i>1)$:  $P_i$ is not identically equal to zero on each $(n-i+1)$-dimensional irreducible  component of the variety defined by a system  $P_1=\dots=P_{i-1}=0$.
\enddefinition

Let  $Y$ be a variety defined in $(\Bbb C^*)^n$ by the system $P_1=\dots=P_k=0$ where $P_1,\dots,P_k$ is a regular sequence. {\it It is easy to see that $\dim Y=n-k$.}

Consider  a variety $X\subset (\Bbb C^*)^n$ defined by an ideal $I$ in the ring $\Cal L_n$ with a basis $Q_1,\dots,Q_N $. Let $\Cal L$ be a space of $\Bbb C$-linear combinations of the functions $Q_1,\dots,Q_N$.
We will prove the following lemma.

\proclaim {Lemma 5} If $P_1,\dots,P_k$ is a generic $k$-tuple of Laurent polynomials belonging to $\Cal L$ and $k=n-\dim X$, then $P_1,\dots,P_k$ is a regular sequence and a variety $Y$  defined by the system $P_1=\dots=P_k=0$
contains the variety $X$.

\endproclaim

 To prove Lemma 5 we will  need the following auxiliary statement.

\proclaim{Lemma 6} Let $A\subset (\Bbb C^*)^n$ be a finite set. Assume that for any $a_i\in A$ there is $P\in \Cal L$ such that $P(a_i)\neq 0$.   Then there is an algebraic hypersurface $\Gamma\subset \Cal L$ such that for any  $Q\in \Cal L\setminus \Gamma$ and for any $a_j\in A$ the inequality $Q(a_j)\neq 0$ holds.
\endproclaim
\demo{Proof} Let $\Gamma_{a_ i}$  be a hyperplane in $\Cal L$ defined by the following condition: $P\in \Gamma_{a_ i}$ if  $P(a_i)=0$. The union $\Gamma= \bigcup_{a_i\in A}\Gamma_{a_ i}$ is an algebraic hypersurface in $\Cal L$. If $Q\notin \Gamma$ then $Q(a_j)\neq 0$ for all $a_j\in A$.
\enddemo

\demo{Proof of Lemma 5} Since $\Cal L\subset I$ all functions $P_1,\dots,P_k$ vanish on $X$ thus $X\subset Y$. We have to explain why a generic sequence is regular. If $X$ coincides with $(\Bbb C^*)^n$ then there is nothing to prove. Otherwise as a first member $P_1$ of a sequence one can take any nonzero element in $\Cal L$. Assume that for  $i<n-k$ we already chose members $P_1, \dots, P_i\in \Cal L$ such that the sequence $P_1,\dots,P_i$ is regular. Consider a variety $Y_i$ defined by the system $P_1=\dots=P_i=0$. The variety $X$ can not contain any irreducible $(n-i)$-dimensional component of $Y_i$ since $\dim X=n-k<n-i$. Take a finite set $A$ containing a point at each $(n-i)$-dimensional component of $Y_i$ not belonging to $X$. According to Lemma 6 there is a hypersurface $\Gamma\subset \Cal L$ such that any $P\in \Cal L \setminus \Gamma$ does not vanish at any point from $A$. As the next member of the sequence one can take any $P_{i+1}\in \Cal L\setminus \Gamma$.
\enddemo

\subhead {4. Developed systems of $n$  equations in $(\Bbb C^*)^n$}
\endsubhead
 Our proof of the good compactification theorem heavily uses developed systems of  $k$ equations in $(\Bbb C^*)^n$. In Newton polyhedra theory there are many results on such systems with $k=n$. Here we recall  these results. This section can be skipped without compromising understanding of the paper.

 Among all systems of $n$ equations $P_1=\dots=P_n=0$ in $(\Bbb C^*)^n$ their is an interesting subclass of {\it developed} systems which  resemble one polynomial equation in one variable. Below we  recall the definition and main properties of such systems.

Let us start with general definitions.  For a a convex polyhedron $\Delta\subset \Bbb R^n$ and a covector $\xi\in (\Bbb R^n)^*$ we  denote by $\Delta^\xi$ the face of $\Delta$ at which the restriction on $\Delta$ of the linear function $ \langle \xi,x\rangle$ attaints its minima. The face $\Delta^\xi$ of the Minkowskii sum $\Delta=\Delta_1+\dots+\Delta_k$ of $k$-tuple $\Delta_1,\dots,\Delta_k$ in $\Bbb R^n$ is equal to  $\Delta^\xi_1+\dots+ \Delta^\xi_k$.

 For a Laurent polynomial $P=\sum a_mx^m$ with Newton polyhedron $\Delta(P)$ we denote by $P^\xi$ the {\it reduction of $P$ in the co-direction $\xi$} defined by the following formula: $P^\xi=\sum _{m\in \Delta^\xi(P)}a_mx^m$. With a system of equations $P_1=\dots=P_k=0$ in $(\Bbb C^*)^n$ and a covector $\xi$ one associates  the reduced in the co-direction $\xi$ system $P_1^\xi=\dots=P_k^\xi=0$.

\definition{Definition 2 {\rm (see [14])}} An $n$-tuple $\Delta_1,\dots\Delta_n$ of convex polyhedra in $\Bbb R^n$ is  {\it developed} if for any nonzero co-vector $\xi\in (\Bbb R^n)^*$ the $n$-tuple $\Delta_1^\xi,\dots\Delta_n^\xi$ contains at least one face $\Delta_j^\xi$ which is a vertex of $\Delta_j$. The system of  equations $P_1 = · · · = P_n = 0 $ in $(\Bbb C^*)^n$
is called {\it developed} if  $n$-tuple $\Delta(P_1,)\dots\Delta(P_n)$ of  Newton polyhedra of $P_1,\dots,P_n$ is developed.
\enddefinition

 A polynomial in one variable
of degree $d$ has exactly d roots counting with multiplicity. The number
of roots in $(\Bbb C^*)^n$
counting with multiplicities of a developed system
is determined only by  Newton polyhedra of equations by the  Bernstein--Koushnirenko formula (if the
system is not developed this formula holds only for generic systems
with given Newton polyhedra).

As in the one-dimensional case, one can explicitly compute the sum
of values of any Laurent polynomial over the roots of a developed system
[6]. In particular,  it allows to eliminate all unknowns but one from the system.

As in the one-dimensional case,  one can explicitly compute the product of all of the roots of the system
regarded as elements in the group $(\Bbb C^*)^n$ [14].

For two polynomials in one variable the following identity holds: up to the sign depending on degrees of the polynomials the product of values of the first polynomial over the roots of the second one is equal  to the product of values of the second polynomial over the roots of the first one multiplied by a certain monomial in coefficients of the polynomials. Assume that for given $(n+1)$ Laurent polynomials $P_1,\dots,P_{n+1}$ in $n$ variable  any $n$-tuple out of $(n+1)$ - tuple of their Newton polyhedra is developed. Then for any $1\leq i<j\leq n+1$ up to sign depending on Newton polyhedra the product of values of $P_i$ over the common roots of all  Laurent polynomials but $P_i$ is equal  to the product of values of $P_j$ over the roots of all  Laurent polynomial but $P_j$ multiplied by a certain monomial in coefficients of all Laurent polynomials. This result and a review of the results mentioned above can be found in  [15].

\subhead {5. Convenient compactifications of $(\Bbb C^*)^n$ for a system of equations}
\endsubhead
With any given system of equations in $(\Bbb C^*)^n$ one can associate a natural class of compactifications of $(\Bbb C^*)^n$ which are  {\it convenient for the system}. In section 10 we construct a good compactification for a variety $X$ as a convenient compactification for some explicitly presented system of equations. In this section we talk about convenient compactifications.

Consider a variety $Y\subset (\Bbb C^*)^n$ defined by a finite system of equations $$P_1=\dots=P_k=0, \tag 1$$ where $P_1,\dots,P_k\in \Cal L_n$. With each convex  polyhedron $\Delta\subset \Bbb R^n$ one can associate its {\it support function} $H_\Delta$ on $(\Bbb R^n)^*$ defined by the relation
$$H_{\Delta}(\xi)=\min_{x\in \Delta}\langle \xi,  x\rangle.$$

\definition {Definition 3 {\rm (see [12--13])}}
 A toric compactification $ M\supset (\Bbb C^*)^n$ is  {\it  convenient  for the system (1)} if support functions  of Newton polyhedra $\Delta_1,\dots,\Delta_k$  of  $P_1,\dots,P_k$ are linear  at each cone $\sigma$ belonging to the fan $\Cal F_M$ of the toric variety $M$.
\enddefinition

It is easy to verify that {\it $M$ is convenient for the system (1) if and only if its fan $\Cal F_M$ is a subdivision of a dual  fan $\Delta^\bot$ to  $\Delta=\Delta(P_1)+\dots+\Delta(P_k)$.}

For some systems  each convenient compactification  provides a good compactification for the variety $Y$ defined by the system. We will show below that a good compactification for any variety $Y\subset (\Bbb C^*)^n$ can be obtained as a convenient compactification for some auxiliary system.

Consider a hypersurface $\Gamma$ in $\Bbb (C^*)^n$ defined by the equation $P=0,$ where $P\in \Cal L_n$.

\proclaim {Lemma 7 {\rm [12--13]}} Any convenient compactification $M$ for the equation $P=0$ provides a good
compactification for  the hypersurface $\Gamma$.
\endproclaim
\demo  {Proof} Let us  show that the closure of $\Gamma$ in $M$ does not contain any  null-orbit. Let $O$ be a null-orbit and let $\sigma$ be the cone in the fan of $M$ corresponding to the affine toric subvariety $M_O$  containing $O$. Since $M$ is a convenient compactification, on the cone $\sigma$ the support function of $\Delta (P)$ is a linear function $\langle \xi, A \rangle$, where $A$ is a vertex of $\Delta (P)$. Let  $\chi_A$ be the character (the monomial) corresponding to the point $A$.  The support function of  $\Delta (P\cdot \chi^{-1}_A)$ equals to zero on $\sigma$ thus $\Delta(P\cdot \chi^{-1})$ belongs to the cone $\sigma^\bot$ dual to $\sigma$, i.e. $P\cdot \chi^{-1}_A$  is regular on $M_O$. Assume that the monomial $\chi_A$ appears in $P$ with the coefficient $C_A$, which is not equal to zero since $A$ is a vertex of $\Delta(P)$.
 The closure of $\Gamma\subset (\Bbb C^*)^n$ in $M_O$ can be defined by the equation $ P\cdot \chi^{-1}_A=0$. It does not contain $O$ since $ P\cdot \chi^{-1}_A(O)=C_A\neq 0.$
\enddemo

The converse  statement to Lemma 7 also is true:
\proclaim{Lemma 8} If $M\supset (\Bbb C^*)^n$ is a good compactification for $\Gamma$ then $M$ is a convenient for the equation $P=0$.
\endproclaim
We will not use Lemma 8 and will not prove it.

\definition {Definition 4 {\rm (see [12--13])}} A system (1) is call {\it $\Delta$-non-degenerate}  if for any covector $\xi\in (\Bbb R^*)^n$ the following condition $(\xi)$ is satisfied:

for  any root $a\in (\Bbb C^*)^n$ of the system $ P_1^\xi=\dots=P_k^\xi=0$ the differentials $d P_1^\xi,\dots,d P_k^\xi$ are independent at the tangent space to $(\Bbb C^*)^n$ at the point $a$.
\enddefinition

\proclaim {Lemma 9 {\rm [12--13]}} Any convenient compactification $M$ for a $\Delta$-non-degenerate system  provides a good compactification for  the variety $Y\subset (\Bbb C^*)^n$ defined by this system.
\endproclaim

If   in the assumptions of Lemma 9 the toric compactification $M$ is smooth then the closure of $Y$ in $M$ is also smooth and transversal to all orbits of $M$ (see [12--13]). One can show that {\it a generic system (1) with the fixed Newton polyhedra $\Delta_1=\Delta(P_1),\dots, \Delta_k=\Delta(P_k)$ is $\Delta$-non-degenerate}. These statements play the key role in  Newton polyhedra theory,  which computes discreet  invariants of  the variety $Y\subset (\Bbb C^*)^n$ in terms of $\Delta_1,\dots,\Delta_k$, where $Y$ is defined by a generic system of equations with  Newton polyhedra $\Delta_1,\dots,\Delta_k$.

Instead of the $\Delta$-non-degeneracy assumption in Lemma 9  one can assume that the ($\xi$) condition (see definition 4) holds only for covectors $\xi$ belonging to cones  in the fan $\Delta^\bot$ for $\Delta=\Delta_1+\dots+\Delta_k$ dual to faces of $\Delta$ whose dimension is smaller then $k$. For $k=1$ this claim coincides with Lemma 7.

In any ideal $I\subset \Cal L_n$ there exists  {\it an universal Grobner basic $P_1,\dots,P_k\in I$} (see for example [17]), a related material can be found in [9]. One can prove the following

\proclaim{Lemma 10} Any convenient compactification $M$ for a system $P_1=\dots=P_k=0,$ where $P_1,\dots,P_k$ is an universal  Grobner basis of an ideal $I\subset \Cal L_n$ provides a good compactification for  the variety $Y\subset (\Bbb C^*)^n$ defined by the ideal $I$.
\endproclaim

Lemma 10 is applicable for any subvariety $Y\subset (\Bbb C^*)^n$ and it provides a standard proof of the good compactification theorem. Usually a universal Grobner basis of an ideal contains a large number of element (see [9]). The universal Grobner basis technique   is far from being  transparent.

 \subhead {6. Developed systems of $k<n$  equations in $(\Bbb C^*)^n$}
\endsubhead
In this section we will deal with  {\it developed   systems} and with complete intersections $Y\subset (\Bbb C^*)^n$ defined by these systems.  A convenient  compactification for a developed system   provides a good compactification for the corresponding complete intersection $Y$.

Let $\Delta_1,\dots,\Delta_k$ be a $k$-tuple of convex polyhedra in $\Bbb R^n$ and let $\Delta$ be the Minkowskii sum $\Delta_1+\dots+\Delta_k$. Each face $\Gamma$ of $\Delta$ is representable as the sum $\Gamma_1+\dots+\Gamma_k$ where $\Gamma_1,\dots,\Gamma_k$ is the (unique) $k$-tuple of faces of these polyhedra.

\definition{Definition 5} A $k$-tuple $\Delta_1,\dots\Delta_k$ of convex polyhedra in $\Bbb R^n$ is  {\it developed} if for any face $\Gamma$ of $\Delta=\Delta_1+\dots+\Delta_k$  such that $\dim \Gamma<k$ in the representation $\Gamma=\Gamma_1+\dots+\Gamma_k$, where $\Gamma_1,\dots,\Gamma_k$ is a   $k$-tuple of faces of $\Delta_1,\dots, \Delta_k$  at least one face $\Gamma_j$  is a vertex of $\Delta_j$. The system (1)
is called {\it developed} if  $k$-tuple $\Delta(P_1,)\dots\Delta(P_n)$ of  Newton polyhedra of $P_1,\dots,P_k$ is developed.
\enddefinition

For $k=n$ Definition 5 is equivalent to Definition 1.

\proclaim{Lemma 11} Let  $Y\subset (\Bbb C^*)^n$ be a complete intersection  defined by a developed system containing $k$ equations and let $M$ be a convenient  compactification for the system. Then the closure of $Y$ in $M$ does not intersect any orbit $O\subset M$ whose dimension is smaller then $k$.
\endproclaim
\demo {Proof} Lemma 11 can be proved in the same way as Lemma 7. Let us show that the closure of $Y$ in $M$ does not intersect any orbits  $O$ in $M$ such that $\dim O <k$.  Let $O$ be  an orbit with  $\dim O<k$ and let $\sigma$ be the cone of dimension $n-\dim O$ in the fan of $M$ corresponding to the smallest affine toric subvariety $M_O$  containing $O$.

 Since $M$ is a convenient compactification the cone   $\sigma$ has to belong to a some cone $\tau$ of the dual fan $\Delta^\bot$ for $\Delta=\Delta_1+\dots+\Delta_k$.

Since $\dim \tau\geq \dim \sigma$  the cone $\tau$ is dual to a face $\Gamma$ of $\Delta$ such that $\dim \Gamma=n-\dim \tau <k$. By assumption in the representation $\Gamma=\Gamma_1+\dots+\Gamma_k$, where each $\Gamma_i$ is a face of $\Delta_i$, there is a face $\Gamma_j$ which is a vertex $A$ of $\Delta_j$. Thus the support function of $\Delta (P_j)$ is the linear function  $\langle \xi, A\rangle$. The vertex $A$ corresponds to the character  $\chi_A$. Assume that the monomial $\chi_A$ appears in $P$ with the coefficient $C_A$, which is not equal to zero since $A$ is a vertex of $\Delta(P)$.
 The closure of $\{P_j=0\} \subset (\Bbb C^*)^n$ in $M_O$ can be defined by the equation $ P_J\cdot \chi^{-1}_A=0$. It does not intersect the orbit $O$ since the function $ P_j\cdot \chi^{-1}_A$ is equal to the constant $C_A\neq 0$ on $O$. The variety $Y$ is contained in the hypersurface $\{P_j=0\}\subset \Bbb C^*)^n$ and can not intersect $O$ as well.

\enddemo

\proclaim{Corollary 12} If  $Y$ is defined by a developed system, then a convenient  compactification for the system is good compactification for $Y$.
\endproclaim
\demo{Proof} Indeed if $Y$ is defined by a system of $k$ equations in the $n$-dimensional torus, then dimension of $Y$ is greater or equal to $n-k$. Thus Corollary 12 follows from Lemma 11.
\enddemo

\subhead {7. Polyhedra with affine independent edges}
\endsubhead
In this section we prove Theorems 15 on geometry  of Newton polyhedra
which easily implies the good compactification theorem.
\definition{Definition 6} A $k$-tuple of segments $I_1,\dots,I_k$ in $\Bbb R^n$ is {\it affine independent} if there is no $k$-tuple of vectors $a_1,\dots, a_k\in \Bbb R^n$ such that the segments $I_1+a_1,\dots,I_k+a_k$ belong to a $(k-1)$-dimensional subspace $L$ of $\Bbb R^n$. Equivalently, $I_1,\dots,I_k$ are {\it affine independent} if dimension of their Minkowski sum $I_1+\dots+I_k$  equals to $k$.
\enddefinition

\definition{Definition 7} Convex polyhedra $\Delta_1,\dots,\Delta_k$ in $\Bbb R^n$ have {\it affine independent edges} if any collection  $I_1\subset \Delta_1,\dots,I_k\subset \Delta_k$ of their edges is affine independent
\enddefinition

\proclaim
{Lemma 13}
Any $k$-tuple of convex polyhedra $\Delta_1,\dots,\Delta_k\subset\Bbb R^n$ having affine independent edges is developed.
\endproclaim
\demo {Proof} Consider the Minkowski sum $\Delta=\Delta_1+\dots+\Delta_k$. Each face $\Gamma$ of $\Delta$ is the Minkowski sum $\Gamma_1+\dots+\Gamma_k$ of the unique collection of faces $\Gamma_1\subset \Delta_1, +\dots, \Gamma_k\subset \Delta_k$ of these polyhedra. If a face $\Gamma_i$ is not a vertex of $\Delta_i,$ then $\Gamma_i$ has to contain an edge $I_i$ of $\Delta_i$. If all faces $\Gamma_1,\dots,\Gamma_k$ are not vertices then $\Gamma_1+\dots+\Gamma_k=\Gamma$ has to contain a sum  $I_1+\dots+I_k$ where $I_1,\dots,I_k$ are some edges of $\Delta_1,\dots,\Delta_k$. Dimension of $I_1+\dots+I_k$ equals to $k$ since $\Delta_1,\dots,\Delta_k$ have affine independent edges. Thus if $\dim \Gamma<k,$ then among faces $\Gamma_1,\dots,\Gamma_k$ has to be at least one vertex.
\enddemo

With any hyperplane $L\subset M$ one can associate a linear function $f_L:M \rightarrow \Bbb R$ (defined up to a nonzero factor) which vanishes on $L$.

\definition {Definition 8} A hyperplane  $L\subset M$ is {\it weakly generic} for a convex polyhedron $\Delta\subset M$ if the restriction of a linear function $f_L$ on  $\Delta$ attaints   its maximum and minimum only at vertices of $\Delta$.
\enddefinition

 Assume that $X\subset (\Bbb C^*)^n$ is   defined by $N+1\geq 1$ equations $$T_1=\dots=T_{N+1}=0,\tag 2$$ (where $T_1,\dots,T_{N+1}$ are Laurent polynomials  not identically equal to zero).

\proclaim {Lemma 14} If a rational hyperplane $L\subset M$ is  weakly generic for the Newton polyhedron $\Delta_1$ of  $T_1,$ then the image $\pi(X)$ of $X$ under the factorization map $\pi: (\Bbb C^*)^n\rightarrow (\Bbb C^*)^n/H(L)=F(L)$ is  an affine subvariety in the torus $F(L)$.  Moreover the dimension of $\pi(X)$ is equal to the dimension of $X$.
\endproclaim
\demo {Sketch of  proof} Conditions on the equation $T_1=0$ and on the hyperplane $L$ imply that the restriction of  $\pi$  on $X$ is a proper map and that each point in $\pi(X)$ has finitely  many   pre-images in $X$. Both claims of Lemma 14
follows from these properties of $\pi$.  We omit  details since in the section 9 we present a constructive proof of Lemma 14.
\enddemo

\definition {Definition 9} A  hyperplane  $L\subset M$ is {\it  generic} for a convex polyhedron $\Delta\subset M$ if the restriction of a linear function $f_L$ on  $\Delta$  is not a constant on any edge of $\Delta$.
\enddefinition

Any hyperplane $L$ generic for $\Delta$  is weakly generic for $\Delta$.

\proclaim {Theorem 15} If  codimension of a subvariety $X$  in the torus $(\Bbb C^*)^n$ is equal to  $k$  then  one can choose a $k$-tuple of Laurent polynomials $P_1,\dots,P_k$ vanishing on $X$ such that their Newton polyhedra $\Delta_1,\dots,\Delta_k$ have affine independent edges.
\endproclaim
\demo{Proof} Induction in codimension  of $X$. If  codimension  is one then $X$ is contained in a hypersurface  $P=0$ where $P$ is a nonzero Laurent polynomial. In that case  one can choose $P_1$ equal to $P$.

If codimension of $X$ is $k>1$ then as $P_1$ one can choose  any nonzero Laurent polynomial vanishing on $X$. After that one can choose a rational hyperplane $L\subset M$ which is  generic for the Newton polyhedron $\Delta_1$ of $P_1$. In particular $L$ is  weakly generic for $\Delta_1$. Thus by Lemma 14 the image $\pi(X)$  of  $X$ is an affine subvariety of codimension $k-1$ in the $(n-1)$-dimensional torus $F(L)$.

The map $\pi^*$ identifies regular functions on $F(L)$ vanishing on $\pi(X)$ with Laurent polynomials on $(\Bbb C^*)^n$  vanishing on $X$ whose Newton polyhedra belong to $L$. By induction there are Laurent polynomials $P_2,\dots,P_k$ vanishing on $X$ whose Newton polyhedra belong to $L$ and  have affine independent edges.   All edges of $\Delta_1$ are not parallel to the space $L$ since  $L$ is generic for $\Delta_1$. Thus the Laurent polynomials $P_1, P_2, \dots, P_k$  have affine independent edges and vanish on $X$. Theorem 15 is proven.

\enddemo
Assume that an algebraic variety $X\subset (\Bbb C^*)^n$ is defined by an ideal $I\subset \Cal L_n$.
\proclaim {Corollary 16 } If  codimension of  $X$  in  $(\Bbb C^*)^n$ is equal to  $k,$  then  one can choose a $k$-tuple of Laurent polynomials $P_1,\dots,P_k$ in the ideal $I$ such that their Newton polyhedra $\Delta_1,\dots,\Delta_k$ have affine independent edges.
\endproclaim
\demo{Proof} According to  Theorem 15 one can chose Laurent polynomials $P_1,\dots,P_k$ vanishing on $X$ such that their Newton polyhedra $\Delta_1,\dots,\Delta_k$ have affine independent edges. By Hilbert's theorem there is a natural number $m$ such that Laurent polynomials $P_1^m,\dots,P_k^m$ belong to the ideal $I$. Newton polyhedra  $m\Delta_1,\dots,m\Delta_k$ of these polynomials have affine independent edges.
\enddemo

The following corollary   provides a version of    the good compactification theorem.

\proclaim {Corollary 17 {\rm (a version of good compactication theorem)}} If  codimension of $X$  in  $(\Bbb C^*)^n$ is equal to  $k$  then  one can choose Laurent polynomials $P_1,\dots,P_k$ vanishing on $X$ such that any convenient compactification for the system  $P_1=\dots=P_k=0$  is a good compactification  for $X$. Moreover for any ideal $I\subset \Cal L_n$ defining  $X$ one can choose  such  Laurent polynomials from the ideal $I$.
\endproclaim
\demo{Proof} By Theorem 15 one can choose Laurent polynomials $P_1,\dots,P_k$ vanishing on $X$ such that their Newton polyhedra $\Delta_1,\dots,\Delta_k$ have affine independent edges. By Lemma 13 the system $P_1=\dots=P_k=0$ is developed. By Lemma 11 any convenient compactification $M$ for the system is a good compactification for the variety defined by the system. By Lemma 4 $M$ is also a good compactification for $X$. By Corollary 15  $P_1,\dots,P_k$ can be chosen in any ideal $I$ defining  $X$.
\enddemo

If  $X \subset (\Bbb C^*)^n$ is defined by an explicitly written system of equations then Theorem 15 and elimination theory allow to construct  a good compactification for $X$ explicitly (see section 10).

\subhead{8. Resultant and elimination of variables}
\endsubhead
In this section we recall classical results on elimination of variables. We adopt  these results for  projections of subvarieties in torus $(\Bbb C^*)^{n} =(\Bbb C^*)^{n-1} \times C^*$  on the first factor $(\Bbb C^*)^{n-1}$.

Let $P=a_0+a_1t+\dots+a_pt^p$ and $Q=b_0+b_1t+\dots+b_qt^q$  be polynomials in $t$ of degrees $\leq p$ and $\leq q$ correspondingly. The resultant $\Cal R_{p,q}(P,Q)$ is a polynomial in $a_0,\dots,a_p,b_0,\dots,b_q$ with integral coefficients. By definition  $\Cal R_{p,q}(P,Q)$ is the determinant of the homogeneous system of  linear equations whose unknowns are the undetermined coefficients of polynomials $\tilde P$ and $\tilde Q$ of degrees $\leq p-1$ and $\leq q-1$ correspondingly satisfying the relation
$$
P\tilde Q=Q\tilde P. \tag 3
$$
Since an ordering of equations in the system is not fixed its determinant $\Cal R_{p,q}(P,Q)$ is defined up to sign.
\proclaim {Lemma 18} If the leading coefficient $a_p$ of the polynomial $P$   (coefficient $b_q$ of the polynomial $Q$) is not equal to zero, then polynomials $P$ and $Q$ have a common factor if and only if their resultant $\Cal R_{p,q}(P,Q)$ is equal to zero.
\endproclaim
\remark {Remark} If $a_p=b_q=0,$ then the resultant $\Cal R_{p,q}(P,Q)$ is equal to zero (even if polynomials $P$ and $Q$ have no common factor).
\endremark

\demo {Proof of Lemma 18} If $P$ and $Q$ have a common factor $T$ with $\deg T>0$ and $P= P_1 T$, $Q= Q_1T$ then the system (3) has a nontrivial solution: one can put $\tilde P=P_1$ and $\tilde Q=Q_1$. Thus the determinant $\Cal R_{p,q}(P,Q)$ equals to zero. On the other hand if the system (3) has a nontrivial solution then $P$ divides $Q\tilde P$. Since $\deg \tilde P<\deg P$ it can happen only if $Q$ and $P$ have a common factor.
\enddemo
Consider a polynomial $P=a_0+a_1t+\dots+a_pt^p$ together with a finite collection of polynomials $Q_i=b_0^i+b_1^it+\dots+b_{q_i}^it^{q_i}$ where $i=1,\dots,N$. Let $Q_\lambda=\lambda_1Q_1+\dots+\lambda_NQ_N$ be a linear combination of polynomials $Q_i$ with  coefficients $\lambda_i$.
\proclaim{Lemma 19} Assume that the leading coefficient $a_p$ of  $P$ is not equal to zero. Let $q=\max _{1\leq i\leq N} q_i$.  Then the polynomials $P, Q_1,\dots, Q_N$ have a common complex root $t_0$ if and only the resultant $\Cal R(\lambda)=\Cal R_{p,q}(P, Q_\lambda)$ is identically equal to zero as a function in $\lambda=(\lambda_1,\dots,\lambda_N)$. If in addition the constant term $a_0$ of $P$ is not equal to zero then any common root of $P, Q_1,\dots Q_N$ also is not equal to zero.
\endproclaim
\demo{Proof} If all polynomials $P, Q_1,\dots, Q_N$ have a common root $t_0$ then for any $N$-tuple $\lambda=(\lambda_1,\dots,\lambda_N)$ the polynomial $Q_\lambda=\lambda_1Q_1+\dots+\lambda_NQ_N$ and the polynomial $P$ have the common root $t_0$. Thus by Lemma 18 the resultant $\Cal R_{p,q}(P,Q_\lambda)=\Cal R(\lambda)$ is identically qual to zero as a function in $\lambda=(\lambda_1,\dots,\lambda_N)$.

Assume now that each root $t_k$, $1\leq k\leq p$ of the polynomial $P$ is not a root of some polynomial $Q_i$. For each root $t_k$ let $\Gamma_k\subset \Bbb C^N$ be the hyperplane  defined by the equation $\sum_{1\leq j\leq N}\lambda_j Q^j(t_k) =0$. If $\lambda^0=(\lambda_1^0,\dots,\lambda_N^0)\in \Bbb C^N$ do not belong to the union $\cup_{1\leq i\leq p } \Gamma_i\subset \Bbb C^N$ then the polynomial $Q_{\lambda^0}=\lambda_1^0 Q^1+\dots+\lambda_N^0 Q^N$ has no common roots with the polynomial $P$. Thus the resultant $\Cal R_{p,q}(P,tQ_\lambda)$ is not identically equal to zero as function in $\lambda$.

If  $a_0\neq 0$ then zero is not a root of $P$ thus it is not a common root of $P,Q_1,\dots,Q_N$.

\enddemo

\subhead {9. Projection of  $X\subset (\Bbb C^*)^{n}$  on a sub-torus  in $(\Bbb C^*)^{n}$}
\endsubhead
In this section we will present a constructive proof of Lemma 14. It will be used as a step in a constructive proof of the good compactification theorem.

\subsubhead {9.1. Modified problem}
\endsubsubhead
Let us modify a little our problem. The group $F(L)$ is the factor-group of $(\Bbb C^*)^n$ by the normal divisor $H(L)$. Let us  choose any complementary to $H(L)$ subgroup $H([e])$ (i.e. such a subgroup that the identity $H([e])\times H(L)=(\Bbb C^*)^n$ holds) and  consider the image $\pi_1(X)\subset H([e])$ under the projection $\pi_1$ of $(\Bbb C^*)^n$ to the first factor $H([e])$.

\proclaim {Problem 1}
In the assumption of Lemma 14 construct explicitly a system of equation
$$
C_m=0, \tag 4
$$
where $C_m \in \Cal L_L$, which defines $\pi_1(X)\subset H([e])$.
\endproclaim

 Any solution of  Problem 1 automatically provides a system of equation on $F(L)$ which defines $\pi (X)\subset F(L)$.
 Indeed, we identified the ring of regular functions on $F_L$ with the ring $\Cal L_L$. Under this identification the system (4) becomes the system of equations on $F_L$ which defines $\pi(X)$.

\remark {Remark} The system of equations (m) constructed below will strongly depend on a choice
of the complementary subgroup $H([e])$ but it defines the same variety $\pi(X)\subset F(L)$.
\endremark

\subsubhead{9.2. Decomposition of $(\Bbb C^*)^n$ into a direct product}
\endsubsubhead
Two rational subspaces $L_1, L_2\subset M$ are {\it complementary to each other} if the  identity  $$\Lambda(L_1) + \Lambda (L_2)=\Lambda $$ holds where
 $\Lambda(L_1) = \Lambda\cap L_1$ and $\Lambda(L_2) = \Lambda\cap L_2$.
If $L_1,L_2$ are complementary to each other then the  identity  $$H(L_1)\times H (L_2)=(\Bbb C^*)^n $$ holds where
 $H(L_1)$ and $H(L_2)$ are subgroups corresponding to the $L_1$ and $L_2$.
 Later we will be interested in the case when $L_1$ is a  hypersurface and $L_2$ is a line.

 We will say that $e$ is a complementary to a rational hypersurface $L$ is $e$ is an irreducible vector in the lattice $\Lambda$ and the line $[e]$ generated by $e$ is a complementary line for $L$. The vector $e$ and the hyperplane $L$ define the linear function $l(e,L):M\rightarrow \Bbb R$ such that $l(e,L)$ vanishes on $L$ and $l(e,L)(e)=1$.

Let $t:(\Bbb C^*)^n\rightarrow  \Bbb C^*$ be the character $\chi_e$ corresponding to the vector $e$.

Let  $t:(\Bbb C^*)^n\rightarrow  \Bbb C^*$ be the character corresponding to the vector $e_t$. Then $t$ has   the following properties:

1) the set $ \{t^{-1}(1)\}\subset (\Bbb C^*)^n$   is the subgroup $H([e])$ in the torus  $(\Bbb C^*)^n$. We identify  regular functions on $H([e])$ (as well as regular functions on  $F(L)$) with  Laurent polynomials from the ring $\Cal L_L$.

2) the map $t:(\Bbb C^*)^n\rightarrow \Bbb C^*$ restricted to  $H(L)$ provides an isomorphism between $H(L)$ and  $\Bbb C^*$. Thus $H(L)$ is a one-parameter group with the parameter $t$.

Each Laurent polynomial on $(\Bbb C^*)^n=H([e])\times H(L)$ can be considered as a Laurent polynomial in  $t$ whose coefficients belong to the ring $\Cal L_L$. We denote by $\pi_1:(\Bbb C^*)^n\rightarrow H([e])$ the  projection of the product to the first factor.

Below we use notations and assumptions from Lemma 14. The lowest degree $m_i$ in $t$ in monomials appearing in $T_i$ is equal to minimum of the function $l_t$ on the Newton polyhedron $\Delta_i$ of $T_i$. Let us put $P=T_1t^{-m_1}, Q_1=T_2t^{-m_2},\dots, Q_{N}=T_{N+1}t^{-m_{N+1}}$. The system (2) defining $X$ is equivalent to the system
$$
P=Q_1=\dots=Q_N=0,
$$
where $P,Q_1,\dots,Q_N$ are  polynomials in $t$  whose coefficients belong to the ring $\Cal L_L$ and in addition the leading coefficient  and the constant term of the polynomial $P$ are characters with  nonzero coefficients. We denote the degrees of the polynomials $P,Q_1,\dots,Q_N$ by $p,q_1,\dots,q_N$ correspondingly

Let $Q_\lambda=\lambda_1Q_1+\dots+\lambda_NQ_N$ be a linear combination of polynomials $Q_i$ with  coefficients $\lambda_i$ and  let $q=\max _{1\leq i\leq N} q_i$. The resultant  $\Cal R_{p,q}(P, Q_\lambda)$ is a polynomial $R(\lambda)$ in $\lambda=(\lambda_1,\dots,\lambda_N)$, i.e.
$$
R(\lambda)=\sum c_{k_1,\dots,k_N}\lambda_1^{k_1}\dots \lambda_N^{k_N},
$$
whose coefficients $c_{k_1,\dots,k_N}$ are Laurent polynomials from the ring $\Cal L_L$.

\proclaim{Theorem 20 {\rm (solution to Problem 1)}} In the assumptions written above  the image $\pi_1(X)$  is an affine subvariety in $H([e])\subset (\Bbb C^*)^n$  defined by the system
$$c_{k_1,\dots,k_N}=0 ,\tag 5$$
where $c_{k_1,\dots,k_N}\in \Cal L_L$ are all coefficients  of the polynomial $\Cal R(\lambda)$. Moreover, $\dim \pi_1(X)=\dim  X$.
\endproclaim

\demo{Proof} We consider functions from the ring $\Cal L_ L$ as functions on $H([e])$.  By assumption the coefficients $a_0$ and   $a_p$ can not vanish at any point $x\in H([e])$. Thus by Lemma 19 the system $P=Q_1=\dots=Q_N=0$ has a common zero $t_0\in \pi_1^{-1}(x)\sim \Bbb C^*$ above a point $x\in H([e])$ if and only if all coefficients $c_{k_1,\dots,k_N}$ vanish at $x$. Any point $x\in \pi_1(X)$ has $\leq p$ pre-images  in $X$  since the degree of $P$ in $t$ is equal to $p$. Thus $\dim \pi_1(X)=\dim X$.
\enddemo
\proclaim {Corollary 21} The image $\pi(X)\subset F(L)$ can be defined by the  system (5) and $\dim \pi(X)=\dim X$
\endproclaim

\subhead{10. Explicit construction of a good compactification}
\endsubhead
Let $X\subset (\Bbb C^*)^n$ be an  algebraic variety. Assume that a system  of equations (2) 
with the following properties is given:

(1)  a variety $Y\subset (\Bbb C^*)^n$ defined by (2) contains the variety $X$;

(2) $\dim Y=\dim X$.

In this section we present an algorithm which replace  (2) by a new system
 $$ P_1=\dots=P_k=0 \tag 6$$
which in addition to the properties (1), (2) has the following properties:

(3) any convenient compactication for (4)  is a good compactification for $X$.

(4) the number of equations in (4) is equal to the codimention of $X$ in $(\Bbb C^*)^n$.

In particular the algorithm allows to compute the dimension of a variety $X$ defined in $(\Bbb C^*)^n$ by a given system of equations and to construct  a good compactification for it.

\subsubhead {10.1. The algorithm}
\endsubsubhead
As the first equation $P_1=0$ of the new system (6) one can take the equation $T_1=0$ (we assume that $T_1$ is not identically equal to zero). To complete the first step of the algorithm some preparations are needed.

Let us choose any rational hyperplane $L_1\subset M$ generic for the Newton polyhedron $\Delta_1$ of $P_1$. Let us choose any complementary vector $e_1$ to  $L_1$. By Theorem 20 one can explicitly write a system of equations $$c^{(1)}_{k_1,\dots,k_N}=0 $$ where $c^{(1)}_{k_1,\dots,k_N} \in \Cal L_{L_1}$ defining the image $Y_1=\pi_1 (Y)\subset H(e_1)$. The first step of the algorithm is completed.
If  all functions $c^{(1)}_{k_1,\dots,k_N}$ are identically equal to zero, then the algorithm  is completed, codimension  of $X$ is one and as the new system (6) contains one equation  $P_1=0$.

The second step is identical to the first step  applied to the system (6) on the torus $H([e_1])$.
In order to do this step we have to replace:

the torus $(\Bbb C^*)^n$ by the torus $H([e_1])\subset (\Bbb C^*)^n$;

the lattice of characters $\Lambda$ by the lattice $\Lambda(L_1)=\Lambda\cap L_1 \subset \Lambda$;

the space of characters $M$ by the space $L_1\subset M$

the ring of Laurent polynomials $\Cal L_n$  by the ring  $\Cal L_{L_1}$.
\smallskip

The subgroup in the torus $H([e_1])$
corresponding to a rational subspace $L\subset L_1$ we will denote by $H_1(L)$.

Let us proceed  with the second step.
If there is a nonzero Laurent polynomial $c^{(1)}_{k_1^0,\dots,k_N^0}$ then as the second equation $P_2=0$ of the new system (6) one can take the equation $c^{(1)}_{k_1^0,\dots,k_N^0}=0$.

After that one can choose a rational hyperplane $L_2$ in the space $L_1$ generic to the Newton polyhedron $\Delta_2\subset L_1$ of $P_2$ and choose a complementary vector $e_2\in \Lambda(L_1)$ to $L_2$ in the space $L_1$.

Let  $t_2:H([e_1])\rightarrow  \Bbb C^*$ be the character corresponding to the vector $e_2$. Then $t_2$ has   the following properties:

1) the set $ \{t^{-1}_2(1)\}\subset H([e_1])$   is the subgroup $H_1([e_2])$ in the torus  $H([e_1])$.    We identify  regular functions on $H_1([e_2])$  with  Laurent polynomials from the ring $\Cal L_{ L_2}$.

2) the map $t_2:H([e_1])\rightarrow \Bbb C^*$ restricted to  $H_1(L_2)$ provides an isomorphism between $H_1(L_2)$ and  $\Bbb C^*$. Thus $H_1(L_2)$ is a one-parameter group with the parameter $t_2$.

The torus $H([e_1])$ can be represented as the product $H_1([e_2])\times H_1(L_2)$. We are interested in the projection $\pi_2(Y_1)\subset H_1(L_2)$ where $\pi_2$ is the projection of the product $H_1([e_2])\times H_1(L_2)$ on the first factor and $Y_1=\pi_1(Y)$.

By Theorem 20 one can explicitly write a system of equations $$c^{(2)}_{k_1,\dots,k_N}=0, $$ where $c^{(2)}_{k_1,\dots,k_N} \in \Cal L_{L_2}$ define the image $Y_2=\pi_2 (Y)\subset H_1(e_2)$. The second step of the algorithm is completed.
If  all functions $c^{(2)}_{k_1,\dots,k_N}$ are identically equal to zero then the algorithm  is completed, codimension  of $X$ is two and as the new system (6) contains two equations  $P_1=P_2=0$.

If there is a nonzero Laurent polynomial $c^{(2)}_{k_1^0,\dots,k_N^0}$ then as the third equation $P_3=0$ of the new system (6) one can take the equation $c^{(2)}_{k_1^0,\dots,k_N^0}=0$ and proceed with the third step of the algorithm and so on. After $k$ steps where $k$ is the codimension of $X$  we will explicitly obtain a system $P_1=\dots=P_k=0$ such that $P_1,\dots,P_k$ vanish on $X$ and their Newton polyhedra have affine independent edges. The description of the algorithm  is completed.

\bigskip
\bigskip
\bigskip

\subhead{11. Modification of Problem 1}
\endsubhead
In this section under assumption of Lemma 14 we consider the following modification of Problem 1.
\proclaim {Problem 2}
Assume that codimension  $X\subset (\Bbb C^*)^n$ is $k>1$. How to construct  a sequence of $(k-1)$ functions
$
R_1,\dots, R_{k-1}
$
from the ring $\Cal L_L$ such that: 1) $R_1,\dots, R_{k-1}$ vanish on $\pi(X)$;
2) $R_1,\dots, R_{k-1}$ form a regular sequence on $H([e])$?
\endproclaim

Theorem 20 suggests the following solution of Problem 2. Consider an auxiliary  linear space $\Cal L$ of $\Bbb C$-linear combinations of the functions $c_{k_1,\dots,k_N}$ (see (5)). As $R_1,\dots, R_{k-1}$ one can take any generic $(k-1)$-tuple of functions from the space $\Cal L$. This solution deal with the space $\Cal L$ of large dimension. Here we present a similar solution which does not involve auxiliary spaces of big dimension.

We  use notations from Lemma 14 and Theorem 20. The variety $X\subset (\Bbb C^*)^n$ defined by the system (2). As the first equation $G_1=0$ we  choose the  equation $T_1=0$. After that we choose a rational hyperplane weakly generic to $\Delta (G_1)$ and  a complementary vector $e$  for $L$. Each function from the ring $\Cal L_n$ can be represented as a Laurent polynomial on $t$ (where $t$ is the character corresponding to the vector $e$) whose coefficients belong to the ring $\Cal L_L$.
One can multiply equations from (2) by any degree of  $t$. Below we  assume that: 1) each $T_i$ is a polynomial in $t$; 2) the polynomial $G_1=T_1$ has a nonzero constant term. Let us denote $\deg G_1$ by $p$ and let $q$ be the maximal degrees of the polynomials $T_i$.

Let $\pi:(\Bbb C^*)^n\rightarrow H(L)$ be the projection defined by the choice of $L$ and $e$.

Let $\Cal L$ be the span of polynomials $T_1,\dots,T_{N+1}$.

Let codimension $X$ be equal to $k$.

We already chose the polynomial $G_1$ from the space $\Cal L$. For any $T\in \Cal L$ by $\Cal R_{p,q}(G_1,T)$ we denote the resultant of $G_1$ and $T$, which we consider as a Laurent polynomial from the ring $\Cal L_L$.

\proclaim {Theorem 22 {\rm (solution to Problem 2)}}  For a generic $(k-1)$-tuple   $G_2,\dots, G_{k}\in \Cal L$ the following conditions hold:

 1)  the sequence  $G_1,G_2,\dots,G_k$ is  regular on $(\Bbb C^*)^n$ and all its members $G_i$ vanish on $X$

2) the sequence  $R_2=\Cal R_{p,q}(G_1, G_2),\dots,R_{k}=\Cal R_{p,q}(G_1, G_{k}) $  is  regular on $H(L)$ and all its members $ R_i$ vanish on $\pi (X)$.
\endproclaim

Our proof of Theorem 22 is similar to the proof of lemma 5. We will  need  an auxiliary Lemma 23 stated below.

Consider a variety  $Y=Z\times \Bbb C^m$ where $Z$ is an affine algebraic variety and $\Bbb C^m$ is a standard linear space with coordinates $\lambda_1, \dots, \lambda_m$. Consider a regular function $\Cal R$ on $Y$ which is a polynomial in $\lambda_1, \dots, \lambda_n$ whose coefficients are regular functions on $Z$.
\proclaim{Lemma 23} Let $A\subset Z$ be a finite set. Assume that for any $a_i\in A$ the restriction of $\Cal R$ on $\{ a_i\}\times \Bbb C^m$ is not identically equal to zero. Then for a generic point $\lambda_0 = (\lambda_{0,1}, \dots, \lambda_{0,n})$ for all points $a_j\in A$ the inequality $\Cal R(a_j,\lambda_0)\neq 0$ holds.
\endproclaim
\demo{Proof} Let $\Gamma_{a_ i}\subset \Bbb C^m$ be a hypersurface in $\Bbb C^m$ defined by the equation $\Cal R(a_i,\lambda)=0$. The union $\Gamma= \bigcup_{a_i\in A}\Gamma_{a_ i}$ is a hypersurface in $\Bbb C^m$. If $\lambda_0\notin \Gamma,$ then $\Cal R(a_j,\lambda_0\neq 0)$ for all $a_j\in A$.
\enddemo

\demo{Proof of Theorem 22} Since $G_1,\dots,G_k \in \Cal L$ all functions $G_i$ vanish on $X$. Assume that for  $1\leq i<k$ we already chose members $G_1, \dots, G_i\in \Cal L$ such that: 1) the sequence $G_1,\dots,G_i$ is regular on $(\Bbb C^*)^n$; 2)  the sequence $\Cal R_2,\dots,\Cal R_i$ is regular on $H(L)$. Consider the variety $Y_i$ defined by the system $G_1=\dots=G_i=0$ on  $(\Bbb C^*)^n$ and the variety $Z_i$ defined by the system $\Cal R_2=\dots=\Cal R_i=0$ on $H(L)$.

The variety $X$ cannot contain any irreducible $(n-i)$-dimensional component of $Y_i$ since $\dim X=n-k<n-i$. Take a finite set $A_i$ containing a point at each $(n-i)$-dimensional component of $Y_i$ not belonging to $X$.

By construction the variety $\pi(X)$ is contained in the variety $Z_i$. he variety $\pi(X)$ cannot contain any irreducible $(n-i-1)$-dimensional component of $Z_i$ since $\dim \pi(X)=n-k-1<n-i-1$. Take a finite set $B_i$ containing a point at each $(n-i-1)$-dimensional component of $Zi$ not belonging to $\pi(X)$.

According to Lemma 6 there is a hypersurface $\Gamma_1\subset \Bbb C^N$ such that for any $\lambda=(\lambda_1,\dots,\lambda_N)$ not in $\Gamma_1$ the function $G=\lambda_1T_1+\dots+\lambda_NT_N $  does not vanish at any point from the set $A_i$.

Consider a function $\Cal R (x,\lambda)$  where $x\in H(L)$ and  $\lambda=(\lambda_1,\dots, \lambda_N)$ defined by the formula $\Cal R (x,\lambda)=\Cal R_{p,q}(G_1, \lambda_1T_1+\dots+\lambda_NT_N)$. Here we consider the resultant  as a polynomial in $\lambda$ whose coefficients belong to the ring $\Cal L_L$ of regular functions on $H(L)$. The restriction of  $\Cal R (x,\lambda)$ to the set  $(x, \lambda)$ with fixed $x=b$ is identically equal to zero if and only if $b\in \pi(X)$.

By the choice of the set $B_i$ and by Lemma 23 there is an algebraic hypersurface $\Gamma_2\subset \Bbb C^N$ such that for any $\lambda$ not in  $\Gamma_2$ the function $\Cal R(x,\lambda)$ doesn't equal to zero at any point of the set $B_i$.

As the next member of the sequence one can take any $G_{i+1}=\lambda_1T_1+\dots+\lambda_NT_N$ for any $ \lambda$ not in $\Gamma_1\cup \Gamma_2$.
\enddemo

\subhead{12. Modification of the algorithm}
\endsubhead
Here we discuss  a modification of the algorithm presented in the section 10.  The modified algorithm does not involve auxiliary spaces of large dimensions.. We assume that the codimension $k$ of $X$ in $(\Bbb C^*)^n$ is given. We will also {\bf make arbitrary generic choices in the construction below}.

Let $X\subset (\Bbb C^*)^n$ be an  algebraic variety of codimension $k$. Assume that a system  of equations $$T_1=\dots=T_{N+1} =0.\tag 7$$
is given which has the following properties :

(1)  a variety $Y\subset (\Bbb C^*)^n$ defined by (7) contains the variety $X$;

(2) $\dim Y=\dim X$.

The modified  algorithm  replaces  (7) by a new system
 $$ P_1=\dots=P_k=0\tag8 $$
containing $k$ equations
which in addition to the properties (1), (2) has the following property:

(3) any convenient compactication for (8)  is a good compactification for $X$.

{\bf The first step of the algorithm}. Consider $\Bbb C$-linear space $\Cal L^0$ consisting of $\Bbb C$-linear combinations $\lambda_1T_1+\dots+\lambda_NT_N$ of the Laurent polynomials $T_1,\dots,T_{N+1}$.
As $P_1$ choose any nonzero element of the space $\Cal L^0$. Choose a rational hyperplane $L$ generic for $\Delta (P_1)$. Choose a complementary vector $e\in \Lambda$ to $L$. Consider the projection $\pi:(\Bbb C^*)^n\rightarrow H(L)$. Take a {generic} $(k-1)$-tuple  $G_2^{(1)},\dots,G_k^{(1)}$ of elements from $\Cal L^0$. Add the first member $G_1^{(1)}=P_1$ to the  $(k-1)$-tuple. Now acting as in Theorem 22 from the sequence $G_1^{(1)},G_2^{(1)},\dots,G_k^{(1)}$ using $\Cal L^0$ and $e$ construct the sequence $R^{(1)}_2,\dots,R^{(1)}_k$ of Laurent polynomials from the space $\Cal L_L$. By Theorem 22 all $R^{(1)}_2,\dots, R^{(1)}_k$ form a regular sequence on $H(L)$ and they vanish on the subvariety $Y_1=\pi(Y)$ of  $H(L)$ having codimension $(k-1)$.

{\bf The second step of the algorithm} is identical to the first step applied to the system $R^{(1)}_2=\dots=R^{(1)}_k=0$ on $H(L), $ where $R^{(1)}_2,\dots,R^{(1)}_k$ belong to $\Cal L_L$ and vanish on the $(k-1)$-dimensional variety $Y_1=\pi(Y)\subset H(L)$.

Consider the $\Bbb C$-linear space $\Cal L^1$ consisting of $\Bbb C$-linear combinations $\lambda_1 R^{(1)}_2+\dots+\lambda_k R^{(1)}_k$.
As $P_2$ choose an element $R^{(1)}_2$ of the space $\Cal L^1$. Choose a rational hyperplane $L_1$  in the space $L$ generic for $\Delta (P_2)\subset L$. Choose a complementary vector $e_1\in \Lambda(L)$ to $L_1$. Consider the projection $\pi_1:H(L)\rightarrow H_1(L_1)$. Take a { generic} $(k-2)$-tuple  $G_3^{(2)},\dots,G_k^{(2)}$ of elements from $\Cal L^1$. Add the second member $G_2^{(2)}=R^{(1)}_2$ to the  $(k-2)$-tuple. Now acting as in Theorem 22 from the sequence $G_2^{(2)},G_3^{(2)},\dots,G_k^{(2)}$ using $L_1$ and $e_1$ construct the sequence $R^{(2)}_3,\dots,R^{(2)}_k$ of Laurent polynomials from the space $\Cal L_{L_1}$. By Theorem 22 all $R^{(2)}_3,\dots, R^{(2)}_k$ form a regular sequence on $H(L_1)$ and they vanish on the subvariety $Y_2=\pi_1(Y_1)$ of $H_1(L_1)$  having codimension $(k-2)$.

Proceeding in the same way one cam make steps $3,\dots,k$. After $k$ steps we obtain a system $P_1,\dots,P_k$ which has needed properties.

\bigskip

{\bf REFERENCES}

\bigskip

[1] D. Bernstein, The number of roots of a system of equations, Funkts. Anal.
Prilozhen. 9, (1975), No. 3, 1--4.

[2] D. Bernstein, A. Kushnirenko, A. Khovanskii, Newton polyhedra, "Uspehi Matem.
Nauk", V. 31, N3, 201--202, 1976 (Russian Math. Surveys).

[3] C. De Concini and C. Procesi, Complete symmetric varieties II Intersection
theory., Adv. Stud. Pure Math., 6 (1985), 481--512.

[4] A. Esterov, B. Kazarnovskii, A.G. Khovanskii, Newton polyhedra and tropical geometry. To
appear in Uspekhi Mat. Nauk (Russian Mathematical Surveys).

[5] W. Fulton, B. Sturmfels, Intersection theory on toric varieties, Topology, 6,
No. 2, (1997), 335--353.

[6] O. Gelfond and A. Khovanskii, Toric geometry and Grothendieck residues, Moscow Mathematical Journal, V. 2, no 1, 2002, 99--112.

[7] I. Itenberg, G. Mikhalkin, E. Shustin, Tropical algebraic geometry, (2nd ed.),
Birkhuser, Basel(2009).

[8] B. Kazarnovskii, Truncations of systems of equations, ideals and varieties, Izv.
Math. 63, No. 3, (1999), 535--547.

[9] B. Kazarnovskii and A. Khovanskii, Tropical Notherian property and Grobner
bases, St. Petersburg Math. J. 26 (2015), No. 5, 797--811.

[10] B. Kazarnovskii and A.G.~Khovanskii, Newton polyhedra, tropical geometry and the ring of condition
for $(C^*)^n$. arXiv:1705.04248.

 [11] G. Kempf, F. Knudsen, D. Mumford, B. Saint-Donat, Toroidal
Embedding. I, Lecture Notes in Math., N 339, Springer-Verlag, 1972.

[12] A. Khovanskii, Newton polyhedra, and toroidal varieties, "Functional Analysis and its
applications", V. 11, N4, 56--64, 1977; translation in Funct. Anal. Appl. 11 (1977), no.
4, 289--296 (1978).

[13] A.G. Khovanskii, Newton polyhedra, and the genus of complete intersections. (Russian)
Funktsional. Anal. i Prilozhen. 12 (1978), no. 1, 51--61.

[14] A. Khovanskii, Newton polyhedrons, a new formula for mixed volume, product of roots
of a system of equations, The Arnoldfest, Proceedings of a Conference in Honor of V.I.
Arnold for his Sixtieth Birthday, Fields Institute Communications, Vol. 24, Amer. Math.
Soc., Providence, RI, 1999, 325--364.

[15] A. Khovanskii, L. Monin, The resultant of developed systems of laurent polynomials.
Moscow Mathematical Journal. V. 17, No. 4, 717--740, 2017.

[16] A. Kouchnirenko, Polyedrea de Newton et nombres de Milnor, Inv. Math. 32,
1--31.

[17] D. Maclagan, B. Sturmfels, Introduction to Tropical Geometry, (Graduate
Studies in Mathematics), (2015).

[18] J. Tevelev, Compactifications of subvarieties of tori American Journal of Mathematics 129, No. 4 (2007), 1087--1104.

\end

\end

\bigskip

\head {CONTENT}
\endhead

\bigskip
1. Introduction
\smallskip

2. Newton polyhedra and dimension of algebraic variety.

(Lemmas 1-2,Theorem 3).
\smallskip

3. Good compactification and regular sequences

(Definition 1, Lemmas 4-6).
\smallskip

4. Developed systems of $n$  equations in $(\Bbb C^*)^n$.

(Definition 2).

\smallskip

5. Convenient compactifications of $(\Bbb C^*)^n$ for a system of equations.

(Definition 3, Lemmas 7-10).
\smallskip

6. Developed systems of $k<n$  equations in $(\Bbb C^*)^n$.

(Definition 4, Lemma 11, Corollary 12).

\smallskip
7. Polyhedra with affine independent edges.

(Definition 5, Lemma 13,  Theorem 14, Corollary 15,16).

\smallskip
8. Elimination of variables.

( Lemmas 17,18   Theorem 19).
\smallskip
9. Projection of  $X\subset (\Bbb C^*)^{n}$  on a sub-torus  in $(\Bbb C^*)^{n}$

10. Explicit construction of good compactification.

\smallskip
11. Modification of Problem 1
\smallskip
12. Modification of the algorithm.

\end